\documentclass{amsart}
\usepackage[all]{xy}
\usepackage{color}
\usepackage{amssymb}

\newtheorem{theorem}{Theorem}[section]
\newtheorem{lemma}[theorem]{Lemma}
\newtheorem{proposition}[theorem]{Proposition}
\newtheorem{corollary}[theorem]{Corollary}

\newtheorem{definition}[theorem]{Definition}
\newtheorem{example}[theorem]{Example}

\theoremstyle{remark}
\newtheorem{notation}[theorem]{Notation}

\theoremstyle{remark}
\newtheorem{remark}[theorem]{Remark}

\newcommand{\ko}{\: , \;}
\newcommand{\ca}{{\mathcal A}}

\newcommand{\cm}{{\mathcal M}}

\newcommand{\cb}{{\mathcal B}}

\newcommand{\cc}{{\mathcal C}}

\newcommand{\co}{{\mathcal O}}

\newcommand{\cd}{{\mathcal D}}

\newcommand{\cv}{{\mathcal V}}
\newcommand{\cw}{{\mathcal W}}

\newcommand{\dgcat}{\mathsf{dgcat}}

\newcommand{\dgcatp}{\mathsf{dgcat}_{\geq 0}}

\newcommand{\internalcomment}[1]{}

\begin{document}

\title[Postnikov towers, $k$-invariants and obstruction theory]{Postnikov towers, $k$-invariants and obstruction theory for DG categories}
\author{Gon{\c c}alo Tabuada}
\address{Departamento de Matematica, FCT-UNL, Quinta da Torre, 2829-516 Caparica,~Portugal}


\email{
\begin{minipage}[t]{5cm}
tabuada@fct.unl.pt
\end{minipage}
}

\begin{abstract}
By inspiring ourselves in Drinfeld's DG quotient, we develop Postnikov towers, $k$-invariants and an obstruction theory for dg categories. As an application, we obtain the following `rigidification' theorem: let $\ca$ be a homologically connective dg category and $F_0:\cb \rightarrow \mathsf{H}_0(\ca)$ a dg functor to its homotopy category. If the family $\{\omega_n(F_n)\}_{n\geq 0}$ of obstruction classes vanishes, then a lift $F:\cb \rightarrow \ca$ for $F_0$ exists.
\end{abstract}

\maketitle

\tableofcontents
\section{Introduction}
A differential graded (=dg) category is a category enriched in the category of complexes of modules over some commutative base ring $R$. Dg categories provide a framework for `homological
geometry' and for non-commutative algebraic geometry in the sense of Bondal, 
Drinfeld, Kapranov, Kontsevich, To{\"e}n, Van den Bergh, $\ldots$ \cite{BV} \cite{BK} \cite{Drinfeld} \cite{Chitalk} \cite{ENS} \cite{IHP} \cite{finMotiv} \cite{Toen}. They are considered as (enriched) derived categories of quasi-coherent sheaves on a hypothetical non-commutative space (see Keller's ICM-talk survey~\cite{ICM}).

In \cite{These}, the homotopy theory of dg categories was constructed. This theory has allowed several developments such as: the creation by To{\"e}n of a derived Morita theory~\cite{Toen}; the construction of a category of `non-commutative motives'~\cite{These}; the first conceptual characterization~\cite{Duke} of Quillen-Waldhausen's $K$-theory~\cite{Quillen1}~\cite{Wald} since its definition in the early $70$'s$\ldots$.

\vspace{0.2cm}

In this article, we develop new `ingredients' in this homotopy theory: Postnikov towers, $k$-invariants and an obstruction theory for homologically connective dg categories.

\vspace{0.2cm}

{\bf Homologically connective dg categories:} A dg category $\ca$ is {\em homologically connective} if for all objects $x,y \in \ca$, the homology $R$-modules $H_i(\ca(x,y))$ are zero for $i<0$. Our motivation comes from non-abelian Hodge theory (see~\cite{Simpson}~\cite{Simpson1}~\cite{Toen1}): 
\begin{example}\label{example1}
We can associate to a $\cc^{\infty}$-manifold $\cm$ its (homologically connective) dg category $T_{DR}(\cm)$ of flat vector bundles on $\cm$. For two flat bundles $V$ and $W$, $T_{DR}(\cm)(V,W)$ is the complex of smooth forms on $\cm$ with coefficients on the vector bundle of morphisms from $V$ to $W$. Although the homotopy category $\mathsf{H}_0(T_{DR}(\cm))$ is equivalent, by the Riemann-Hilbert correspondance, to the category of finite dimensional linear representations (up to homotopy) of the fundamental group of $\cm$, the dg category $T_{DR}(\cm)$ carries much more information: for two flat bundles $V$ and $W$, corresponding to two local systems $L_1$ and $L_2$, the homology group $\mathsf{H}_i(T_{DR}(\cm)(V,W))$ is isomorphic to the $\mbox{Ext}$-group $\mbox{Ext}_i(L_1, L_2)$, computed in the category of abelian sheaves over $\cm$.

\vspace{0.1cm}

We can also associate to a complex manifold $\cm$ its (homologically connective) dg category $T_{Dol}(\cm)$ of holomorphic complex vector bundles on $\cm$. For two holomorphic bundles $V$ and $W$, $T_{Dol}(\cm)(V,W)$ is the Dolbeault complex with coefficients in the vector bundle of morphisms from $V$ to $W$. Although the homotopy category $\mathsf{H}_0(T_{Dol}(\cm))$ is equivalent to the full subcategory of the bounded coherent derived category of $\cm$, whose objects are the holomorphic vector bundles, the dg category $T_{Dol}(\cm)$ carries much more information: for two holomorphic vector bundles $V$ and $W$, the homology group $\mathsf{H}_i(T_{Dol}(\cm)(V,W))$ is isomorphic to the $\mbox{Ext}$-group $\mbox{Ext}^i_{\co}(\cv, \cw)$ calculated in the category of holomorphic coherent sheaves. For instance, if ${\bf 1}$ if the trivial vector bundle of rank $1$ and $V$ is any holomorphic vector bundle
$$ \mathsf{H}_i(T_{Dol}(\cm)({\bf 1},V)) \simeq \mathsf{H}_i(X,V)\,,$$
where $\mathsf{H}_i(X,V)$ is the $i$th Dolbeault homology group of $V$.
\end{example}
The purpose of this article is to develop a general `technology' that allow us to characterize precisely which are the obstructions appearing when one tries to lift `information' from the homotopy category to the differential graded one. 
\vspace{0.2cm}

{\bf Postnikov towers:} A {\em Postnikov tower} $(\ca_n)_{n \geq 0}$ for a homologically connective dg category $\ca$ is a commutative diagram in the category $\dgcat$ of dg categories
$$
\xymatrix@!0 @R=2.5pc @C=4pc{
&& \vdots \ar[d] \\
&& \ca_2 \ar[d] \\
&& \ca_1 \ar[d] \\
\ca \ar[rr]^{P_0} \ar[urr]^{P_1} \ar[uurr]^{P_2} & &\ca_0
}
$$
such that:
\begin{itemize}
\item[A)] The dg functor $P_n: \ca \longrightarrow \ca_n$ satisfies the following conditions:
\begin{itemize}
\item[A1)] for all objects $x,y \in \ca$, the induced map on the homology $R$-modules
$$ \mathsf{H}_i(\ca(x,y)) \stackrel{\sim}{\longrightarrow} \mathsf{H}_i(\ca_n(P_nx,P_ny))\,,$$
is an isomorphism for $i \leq n$ and
\item[A2)] the dg functor $P_n$ induces an equivalence of categories $\mathsf{H}_0(\ca) \stackrel{\sim}{\rightarrow} \mathsf{H}_0(\ca_n)$.
\end{itemize}
\item[B)] For all objects $x,y \in \ca_n$, the homology $R$-modules $\mathsf{H}_i(\ca_n(x,y))$ are zero for $i>n$.
\end{itemize}
By inspiring ourselves in Drinfeld's description of the Hom complexes in his DG quotient (see~\cite[3.1]{Drinfeld}), we construct in section \ref{secBig} a Big (functorial) Postnikov model $P(\ca)$ for $\ca$. We then use it to prove the following `uniqueness' theorem.
\begin{theorem}(\ref{uniqueness})
Given two objects in the category $\mbox{Post}(\ca)$ of Postnikov towers for $\ca$, there exists a zig-zag of weak equivalences relating the two.
\end{theorem}
For many purposes, a dg category $\ca$ can be replaced by any of its {\em Postnikov sections} $\ca_n$. For example if one is only interested in its homotopy category $\mathsf{H}_0(\ca)$ or if one is only interested in its homology $R$-modules in a finite range of dimensions. 

On the other hand, using a small Postnikov model $\mathbb{P}(\ca)$ for $\ca$ (see~\ref{smallsec}), we prove that the full homotopy type of $\ca$ can be recovered from any of its Postnikov towers by a homotopy limit procedure (see proposition~\ref{homtype}).

\vspace{0.2cm}

{\bf $k$-invariants:} Having seen how to decompose a homologically connective dg category $\ca$ into its Postnikov sections $\ca_n, \, n \geq 0$, we consider the inverse problem of building a Postnikov tower for $\ca$, starting with $\ca_0$ and inductively constructing $\ca_{n+1}$ from $\ca_n$. In order to solve this problem, we construct (see~\ref{k-invariant}) a dg functor 
$$ \gamma_n: P_n(\ca) \longrightarrow \mathbb{P}_n(\ca)\ltimes\mathsf{H}_{n+1}(\ca)[n+2]\,,$$
from the $n$th Big Postnikov section of $\ca$ to a square zero extension (see~\ref{Eilenberg}) of $\mathbb{P}_n(\ca)$. The image  of $\gamma_n$ in the homotopy category $\mathsf{Ho}(\dgcat \downarrow \mathbb{P}_n(\ca))$ of dg categories over $\mathbb{P}_n(\ca)$ is called the {\em $n$th $k$-invariant $\alpha_n(\ca)$ of $\ca$} (see~\ref{def1}). We show that $\alpha_n(\ca)$ corresponds to a derived derivation of $\mathbb{P}_n(\ca)$ with values in the $\mathbb{P}_n(\ca)\text{-}\mathbb{P}_n(\ca)$-bimodule $\mathsf{H}_{n+1}(\ca)[n+2]$ (see~\ref{Dderivations}).

Then we prove our main theorem, which shows how the full homotopy type of $P_{n+1}(\ca)$ in $\dgcat$ can be entirely recovered from $\alpha_n(\ca)$.
\begin{theorem}(\ref{mainthm})
We have a homotopy fiber sequence
$$ P_{n+1}(\ca) \longrightarrow P_n(\ca) \stackrel{\gamma_n}{\longrightarrow} \mathbb{P}_n(\ca)\ltimes \mathsf{H}_{n+1}(\ca)[n+2]$$
in $\mathsf{Ho}(\dgcat\downarrow \mathbb{P}_n(\ca))$.
\end{theorem}

\vspace{0.2cm}

{\bf Obstruction theory:} By inspiring ourselves in the examples appearing in non-abelian Hodge theory, we formulate the following general `rigidification' problem:

\vspace{0.02cm}

Let $\ca$ be an homologically connective dg category and $F_0: \cb \longrightarrow \mathsf{H}_0(\ca)$ a dg functor with values in its homotopy category, with $\cb$ a cofibrant dg category. Is there a lift $F: \cb \longrightarrow \ca$ making the diagram
$$
\xymatrix@!0 @R=2.5pc @C=4pc{
& \ca \ar[d]^{\tau_{\leq 0}} \\
\cb \ar[r]_-{F_0} \ar[ur]^F & \mathsf{H}_0(\ca)
}
$$
commute?

\vspace{0.1cm}

Intuitively the dg functor $F_0$ represents the `up-to-homotopy' information that one would like to rigidify, i.e. lift to the dg category $\ca$. 

In order to solve this problem, we consider a Postnikov tower for $\ca$ (e.g. its Big Postnikov model)
$$
\xymatrix{
&& \vdots \ar[d]\\
&& P_2(\ca) \ar[d]\\
&& P_1(\ca) \ar[d]\\
\cb \ar@{-->}[uurr]^-{F_2} \ar@{-->}[urr]^-{F_1} \ar[rr]^-{F_0} && \mathsf{H}_0(\ca) \simeq P_0(\ca)
}
$$
and we try to lift $F_0$ to dg functors $F_n:\cb\rightarrow P_n(\ca)$ for $n=1,2 \ldots $ in sucession. The image of the composed dg functor
$$ \cb \stackrel{F_n}{\longrightarrow} P_n(\ca) \stackrel{\gamma_n}{\longrightarrow} \mathbb{P}_n(\ca) \ltimes \mathsf{H}_{n+1}(\ca)[n+2]$$
in the homotopy category $\mathsf{Ho}(\dgcat \downarrow \mathbb{P}_n(\ca))$ is called the {\em obstruction class $\omega_n(F_n)$ of $F_n$} (see~\ref{obstruction}). We interpret it as a derived derivation of $\cb$ with values in a $\cb\text{-}\cb$-bimodule (see~\ref{hoje3}).

We then prove that, if at each stage of the inductive process of constructing lifts $F_n: \cb \longrightarrow \mathbb{P}_n(\ca)$, the obstruction class $\omega_n(F_n)$ vanishes, then a lift $F:\cb \longrightarrow \ca$ for $F_0$ exists.

\begin{theorem}(\ref{final})
If the family $\{\omega_n(F_n)\}_{n \geq 0}$ of obstruction classes vanishes, then the `rigidification' problem has a solution.
\end{theorem}

\section{Acknowledgements}
It is a great pleasure to thank Carlos Simpson for motivating conversations and Gustavo Granja for several important comments on an older version of this article. I would like also to thank the Laboratoire J.-A. Dieudonn{\'e} at Nice-France for his hospitality, where some of this work was carried out. 

\section{Preliminaries}
In what follows, $R$ will denote a commutative ring with unit. The
tensor product $\otimes$ will denote the tensor product over $R$. Let $Ch$ be the category of complexes of $R$-modules and $Ch_{\geq 0}$ the full subcategory of positively graded complexes (we consider homological notation, i.e. the differential decreases the degree). Recall from~\cite[2.3.11]{HoveyLivro}, that $Ch$ carries a projective model structure, whose weak equivalences are the quasi-isomorphisms and whose fibrations are the degreewise surjective maps. 

We denote by $\dgcat$ the category of small dg categories, see \cite{Drinfeld}~\cite{ICM}~\cite{These}.
\begin{definition}\label{bimodule}
Let $\ca$ be a small dg category.
\begin{itemize}
\item[-] the {\em opposite dg category} $\mathcal{A}^{op}$ of $\ca$ has the same objects as $\mathcal{A}$ and its complexes on morphisms are defined by $\mathcal{A}^{op}(x,y) =\mathcal{A}(y,x)$.
\item[-] a {\em $\ca\text{-}\ca$-bimodule $M$} is a dg functor $M:\ca^{op}\otimes\ca \rightarrow Ch$.
\end{itemize}
\end{definition}

 Recall from~\cite[1.8]{These} that $\dgcat$ carries a cofibrantly generated Quillen model structure whose weak equivalences are defined as follows:
\begin{definition}
A dg functor $F:\ca \longrightarrow \cb$ is a {\em quasi-equivalence} if:
\item[(i)] for all objects $x,y \in \ca$, the induced morphism 
$$F(x,y): \ca(x,y) \stackrel{\sim}{\longrightarrow} \cb(Fx,Fy)$$
 is a quasi-isomorphism in $Ch$ and
\item[(ii)] the induced functor $\mathsf{H}_0(F): \mathsf{H}_0(\ca) \stackrel{\sim}{\longrightarrow} \mathsf{H}_0(\cb)$ is an equivalence of categories.
\end{definition}
\begin{remark}\label{remquasi}
Notice that if condition $(i)$ is verified, condition $(ii)$ is equivalent to:
\begin{itemize}
\item[(ii)'] the induced functor
$$\mathsf{H}_0(F): \mathsf{H}_0(\ca) \longrightarrow \mathsf{H}_0(\cb)$$
is essentially surjective.
\end{itemize}
\end{remark}
Let us now recall from \cite[1.13]{These}, the following characterization of the fibrations in $\dgcat$.
\begin{proposition}\label{deffibrations}
A dg functor $F:\ca \longrightarrow \cb$ is a fibration if and only if:
\begin{itemize}
\item[F1)] for all objects $x,y \in \ca$, the induced morphism 
$$\xymatrix{F(x,y): \ca(x,y) \ar@{->>}[r] & \cb(Fx,Fy)}$$
 is a fibration in $Ch$ and
\item[F2)] for every objet $a_1 \in \ca$ and every morphism $v
  \in \cb(F(a_1),b)$ which becomes invertible in $\mathsf{H}_0(\mathcal{B})$, there exists a morphism $u \in
    \ca(a_1,a_2)$ such that $F(u)=v$ and which become invertible in $\mathsf{H}_0(\mathcal{A})$.
\end{itemize}
\end{proposition}
\begin{remark}\label{objfib}
Since the terminal object in $\dgcat$ is the zero category $0$ (one object and trivial dg algebra of endomorphisms), every object in $\dgcat$ is fibrant. 
\end{remark}
\begin{corollary}\label{corfib}
Let $F:\ca \longrightarrow \cb$ be a dg functor such that:
\begin{itemize}
\item[-] induces a surjective map on the set of objects,
\item[-] for all objects $x,y \in \ca$, the induced morphism
$$\xymatrix{ F(x,y): \ca(x,y) \ar@{->>}[r] & \cb(Fx,Fy)}$$
is a fibration in $Ch$ and
\item[-] the induced functor
$$ \mathsf{H}_0(F): \mathsf{H}_0(\ca) \stackrel{\sim}{\longrightarrow} \mathsf{H}_0(\cb)$$
is an equivalence of categories.
\end{itemize}
Then $F$ is a fibration in $\dgcat$.
\end{corollary}
\begin{definition}\label{connected}
Let $\ca$ be a small dg category.
\begin{itemize}
\item[-] We say that $\ca$ is {\em homologically connective} if for all objects $x,y \in \ca$, the homology $R$-modules $H_i(\ca(x,y))$ are zero for $i<0$.
\item[-] We say that $\ca$ is {\em positively graded} if for all objects $x,y \in \ca$, the $R$-modules $\ca(x,y)_i$ are zero for $i<0$.
\end{itemize}
\end{definition}
\begin{notation}
We denote by $\dgcatp$ the category of small positively graded dg categories.
\end{notation}
Recall from~\cite[4.16]{DGScat} that we have an adjunction
$$
\xymatrix{
\dgcat \ar@<1ex>[d]^{\tau_{\geq 0}} \\
\dgcatp \ar@{^{(}->}@<1ex>[u]^i\,,
}
$$
where $\tau_{\geq 0}$ denotes the `intelligent' truncation functor.
\begin{remark}\label{pos->con}
Notice that for a homologically connective dg category $\ca$, the co-unit of the previous adjunction, furnishes us a natural quasi-equivalence
$$ \eta_{\ca}: \tau_{\geq 0}(\ca) \stackrel{\sim}{\longrightarrow} \ca\,,$$
which induces the identity map on set of objects.
This (functorial) procedure will allow us to extended several constructions from positively graded to homologically connective dg categories.
\end{remark}
We finish these preliminaries with some homotopical algebra results and the notion of lax monoidal functor. Let $\cm$ be a Quillen model category and $X$ an object of $\cm$.
\begin{notation}\label{overcat}
We denote by $\cm \downarrow X$ the category of objects of $\cm$ over $X$, see~\cite[7.6.2]{Hirschhorn}. Notice that its terminal object is the identity morphism on $X$. 
\end{notation}
\begin{remark}\label{rks}
Recall from \cite[7.6.5]{Hirschhorn} that $\cm \downarrow X$ carries a natural Quillen model structure induced by the one on $\cm$. In particular an object $Y \longrightarrow X$ in $\cm \downarrow X$ is cofibrant if and only if $Y$ is cofibrant in $\cm$ and is fibrant if and only if the morphism $\xymatrix{Y \ar@{->>}[r] & X}$ is a fibration in $\cm$. 
Notice also that if $f:X \longrightarrow X'$ is a morphism in $\cm$, we have a Quillen adjunction
$$
\xymatrix{
\cm \downarrow X \ar@<-1ex>[d]_{f_!} \\
\cm \downarrow X' \ar@<-1ex>[u]_{f^!}\,,
}
$$
where $f^!$ associates to an object $Y \longrightarrow X'$ in $\cm \downarrow X'$ the object $X \underset{X'}{\times} Y \longrightarrow X$ in $\cm\downarrow X$ and $f_!$ associates to an object $Z \longrightarrow X$ in $\cm\downarrow X$ the object $Z \rightarrow X \stackrel{f}{\longrightarrow} X'$ in $\cm\downarrow X'$.
We have also a natural forgetful functor
$$ U: \cm \downarrow X \longrightarrow \cm\,,$$
which preserves cofibrations, fibrations and weak equivalences. This implies that $U$ descends to the homotopy categories $U: \mathsf{Ho}(\cm \downarrow X) \longrightarrow \mathsf{Ho}(\cm)$ and so we obtain the following lemma.
\end{remark}
\begin{lemma}\label{rks1}
Let $f$ and $f'$ be two morphisms in $\cm \downarrow X$. If they become equal in $\mathsf{Ho}(\cm \downarrow X)$, then $U(f)$ and $U(f')$ become equal in $\mathsf{Ho}(\cm)$.
\end{lemma}

\begin{lemma}\label{strictification}
Let $\cm$ be a Quillen model category. Suppose we have a (non-commutative) diagram
$$
\xymatrix{
& X \ar@{->>}[d]^p \\
Z \ar[ur]^{f'} \ar[r]_f & Y\,,
}
$$
where $Z$ is cofibrant, $Y$ is fibrant, $p$ is a fibration in $\cm$ and the composition $p\circ f'$ becomes equal to $f$ in the homotopy category $\mathsf{Ho}(\cm)$. Then, there exists a lift $\widetilde{f}:Z \longrightarrow X$ of $f$ which makes the diagram
$$
\xymatrix{
& X \ar@{->>}[d]^p \\
Z \ar[ur]^{\widetilde{f}} \ar[r]_f & Y\
}
$$
commute.
\end{lemma}
\begin{proof}
Notice that since $Z$ is cofibrant and $Y$ is fibrant, the composition $p\circ f'$ becomes equal to $f$ in $\mathsf{Ho}(\cm)$, if and only if $p\circ f'$ and $f$ are left homotopic. This allow us to construct a (solid) commutative square
$$
\xymatrix{
*+<1pc>{Z} \ar@{>->}[d]^{\sim}_{i_0} \ar[r]^{f'} & X \ar@{->>}[d]^p \\
I(Z) \ar[r]_H \ar@{-->}[ur]^{\widetilde{H}} & Y\,,
}
$$
where $I(Z)$ is a cylinder object for $Z$ and $H$ is an homotopy between $p \circ f'$ and $f$. Finally, $p$ has the right lifting property with respect to $i_0$ and so we obtain a desired morphism
$$ \widetilde{f}: Z \stackrel{i_1}{\longrightarrow} I(Z) \stackrel{\widetilde{H}}{\longrightarrow} X\,,$$
such that $p \circ \widetilde{f} =f$.
\end{proof}

\begin{definition}\label{lax}
Let $(\cc, -\otimes-, \mathbb{I}_{\cc})$ and $(\cd, -\wedge-, \mathbb{I}_{\cd})$ be two symmetric monoidal categories. A {\em lax monoidal functor} is a functor $F:\cc \longrightarrow \cd$ equipped with:
\begin{itemize}
\item[-] a morphism $\eta: \mathbb{I}_{\cd} \longrightarrow F(\mathbb{I}_{\cc})$ and 
\item[-] natural morphisms
$$ \psi_{X,Y}: F(X) \wedge F(Y) \longrightarrow F(X \otimes Y),\,\,\,\, X,Y \in \cc$$
which are coherently associative and unital (see diagrams $6.27$ and $6.28$ in \cite{Borceaux}).
\end{itemize}
A lax monoidal functor is {\em strong monoidal} if the morphisms $\eta$ and $\psi_{X,Y}$ are isomorphisms.
\end{definition}
Throughout this article the adjunctions are displayed vertically with the left, resp. right, adjoint on the left side, resp. right side.
\section{Postnikov towers}
In this chapter, we construct (functorial) Postnikov towers for homologically connective dg categories. We prove that they are `essentially' unique (see theorem~\ref{uniqueness}) and that the full homotopy type of a homologically connective dg category can be recovered from any of its Postnikov towers (see proposition~\ref{homtype}).
\begin{definition}\label{maindef}
A {\em Postnikov tower} $(\ca_n)_{n \geq 0}$ for a positively graded dg category $\ca$ is a commutative diagram in $\dgcat$
$$
\xymatrix{
& &\vdots \ar[d] \\
& &\ca_2 \ar[d] \\
& &\ca_1 \ar[d] \\
\ca \ar[rr]^{P_0} \ar[urr]^{P_1} \ar[uurr]^{P_2} & &\ca_0
}
$$
such that:
\begin{itemize}
\item[A)] The dg functor $P_n: \ca \longrightarrow \ca_n$ satisfies the following conditions:
\begin{itemize}
\item[A1)] for all objects $x,y \in \ca$, the induced map on the homology $R$-modules
$$ \mathsf{H}_i(\ca(x,y)) \stackrel{\sim}{\longrightarrow} \mathsf{H}_i(\ca_n(P_nx,P_ny))$$
is an isomorphism for $i \leq n$ and
\item[A2)] it induces an equivalence of categories $\mathsf{H}_0(\ca) \stackrel{\sim}{\longrightarrow} \mathsf{H}_0(\ca_n)$.
\end{itemize}
\item[B)] For all objects $x,y \in \ca_n$, the homology $R$-modules $\mathsf{H}_i(\ca_n(x,y))$ are zero for $i>n$.
\end{itemize}
The dg functor $P_n:\ca \longrightarrow \ca_n$ is called the {\em $n$th Postnikov section} of $\ca$. 
\end{definition}
\begin{remark}\label{H0}
By the $2$ out of $3$ property, the dg functors $\ca_{n+1} \longrightarrow~\ca_n$ induce an equivalence of categories $\mathsf{H}_0(\ca_{n+1}) \stackrel{\sim}{\longrightarrow} \mathsf{H}_0(\ca_n)$.
\end{remark}
\begin{definition}
A {\em morphism $M: (\ca_n)_{n \geq 0} \longrightarrow (\ca_n')_{n \geq 0}$} between two Postnikov towers for $\ca$ is a family of dg functors $M_n:\ca_n \longrightarrow \ca_n'$ which makes the obvious diagrams commute.
\end{definition}
\begin{notation}\label{catPost}
We denote by $\mbox{Post}(\ca)$ the category of Postnikov towers for $\ca$. 
\end{notation}
\begin{remark}\label{morph}
Let $M: (\ca_n)_{n \geq 0} \longrightarrow (\ca_n')_{n \geq 0}$ be a morphism between Postnikov towers for $\ca$. By the 2 out of 3 property, its Postnikov sections $M_n: \ca_n \stackrel{\sim}{\longrightarrow} \ca_n'$ are all quasi-equivalences.
\end{remark}
\begin{remark}\label{fibreplace}
Observe that in a Postnikov tower $(\ca_n)_{n\geq 0}$ for $\ca$, we can replace each dg functor $\ca_{n+1} \longrightarrow \ca_{n}$ by a fibration $\xymatrix{F(\ca_{n+1}) \ar@{->>}[r] &F(\ca_{n})}$, starting with $\ca_1 \longrightarrow \ca_0$ and then going upward. For the inductive step, we factor the composition $\xymatrix{\ca_{n+1} \ar[r] & *+<1pc>{\ca_{n}} \ar@{>->}[r]^-{\sim} & F(\ca_{n})}$ by a trivial cofibration followed by a fibration $\xymatrix{ F(\ca_{n+1}) \ar@{->>}[r] & F(\ca_{n})}$. We obtain then a morphism $(\ca_n)_{n\geq 0} \longrightarrow F(\ca_n)_{n \geq 0}$ between Postnikov towers
$$
\xymatrix{
& &\vdots \ar[d] & \vdots \ar@{->>}[d] \\
& &*+<1pc>{\ca_2} \ar@{>->}[r]^-{\sim} \ar[d] & F(\ca_2) \ar@{->>}[d]\\
& &*+<1pc>{\ca_1} \ar[d] \ar@{>->}[r]^-{\sim} & F(\ca_1) \ar@{->>}[d] \\
\ca \ar[rr]^{P_0} \ar[urr]^{P_1} \ar[uurr]^{P_2} && \ca_0 \ar@{=}[r] & \ca_0 \,.
}
$$
\end{remark}

\begin{definition}\label{maindef1}
Let $\ca$ be a homologically connective dg category. By a Postnikov tower for $\ca$, we mean a Postnikov tower for $\tau_{\geq 0}(\ca)$, see remark~\ref{pos->con}.
\end{definition}
We now present two functorial Postnikov tower models.

\subsection{Small model}\label{smallsec}
Let $n \geq 0$. Consider the `intelligent' truncation functor
$$ \tau_{\leq n}: Ch_{\geq 0} \longrightarrow Ch_{\geq 0}$$
which associates to a complex
$$ M_{\bullet}: \,\,\,\,\, 0 \leftarrow M_{0} \leftarrow \cdots \leftarrow M_{n-1} \leftarrow M_n \leftarrow M_{n+1} \leftarrow \cdots $$
its `intelligent' truncation
$$ \tau_{\leq n}(M_{\bullet}): \,\,\,\,\, 0 \leftarrow M_0 \leftarrow \cdots \leftarrow M_{n-1} \leftarrow M_n/Im(M_{n+1}) \leftarrow 0 \leftarrow \cdots\,.$$
Notice that when $n$ varies, we obtain the following natural tower of complexes
$$
\xymatrix{
& &\vdots \ar@{->>}[d] \\
& &\tau_{\leq 2}(M_{\bullet}) \ar@{->>}[d] \\
& & \tau_{\leq 1}(M_{\bullet}) \ar@{->>}[d]  \\
M_{\bullet} \ar[rr] \ar[urr] \ar[uurr] & & \tau_{\leq 0}(M_{\bullet})\,.
}
$$
Moreover each vertical map is a fibration and the induced map on the homology $R$-modules
$$ \mathsf{H}_i(M_{\bullet}) \stackrel{\sim}{\longrightarrow} \mathsf{H}_i(\tau_{\leq n}(M_{\bullet}))$$
is an isomorphism for $i\leq n$. Notice also that the homology $R$-modules $\mathsf{H}_i(\tau_{\leq n}(M_{\bullet}))$ are zero for $i >n$. 

Now, let $\ca$ be a positively graded dg category. Since for every $n\geq 0$, the truncation functor $\tau_{\leq n}$ is lax monoidal (see~\ref{lax}), the above remarks imply the following: if we apply the `intelligent' truncation functors to each complex of morphisms of $\ca$, we obtain a Postnikov tower
$$
\xymatrix{
& &\vdots \ar@{->>}[d] \\
& &\tau_{\leq 2}(\ca) \ar@{->>}[d] \\
& &\tau_{\leq 1}(\ca) \ar@{->>}[d] \\
\ca \ar[rr] \ar[urr] \ar[uurr] & & \tau_{\leq 0}(\ca) 
}
$$
for $\ca$. Moreover, by construction, all the dg functors in the diagram induce the identity map on the set of objects.
Notice also that since the morphisms of complexes
$$\xymatrix{ \tau_{\leq n+1}(M_{\bullet}) \ar@{->>}[r] & \tau_{\leq n}(M_{\bullet})}$$
are fibrations, remark~\ref{H0} and corollary~\ref{corfib}, imply that the dg functors
$$\xymatrix{ \tau_{\leq n+1}(\ca) \ar@{->>}[r] & \tau_{\leq n}(\ca)}$$
are fibrations in $\dgcat$.
\begin{notation}\label{smallnot}
We denote by $\mathbb{P}(\ca)$ the small Postnikov model obtained. In particular $\mathbb{P}_n(\ca)$ denotes the dg category $\tau_{\leq n}(\ca)$.
\end{notation}
\subsection{Big model}\label{secBig}
We start by recalling from \cite[1.3]{These} same generating cofibrations for the Quillen model structure on $\dgcat$. 
\begin{definition}\label{novadef}
For $n \in \mathbb{Z}$, let $S_{n}$ be the complex $R[n]$ (with $R$ concentrated in degree $n$) and let $D_{n+1}$ be the mapping cone on the identity of $S_{n}$. We denote by ${\bf 1}_n$ the element of degree $n$ in $S_n$, which corresponds to the unit of $R$. Let $\mathcal{C}(n)$ be the dg category with two objects $1$ et $2$ such that
$
\mathcal{C}(n)(1,1)=R \ko
\mathcal{C}(n)(2,2)=R \ko
\mathcal{C}(n)(2,1)=0  \ko
\mathcal{C}(n)(1,2)=S_n
$ and composition given by multiplication.
We denote by $\mathcal{D}(n+1)$ the dg category with two objects $3$ and $4$ such that
$
\mathcal{P}(n)(3,3)=R \ko
\mathcal{P}(n)(4,4)=R \ko
\mathcal{P}(n)(4,3)=0 \ko
\mathcal{P}(n)(3,4)=D_{n+1}
$ and with composition given by multiplication. Finally, let $S(n)$ be the dg functor from $\mathcal{C}(n)$ to
$\mathcal{D}(n+1)$ that sends $1$ to $3$, $2$ to $4$ and $S_{n}$ to $D_{n+1}$
by the identity on $R$ in degree $n$.
\end{definition}
\begin{lemma}\label{lift}
Let $\ca$ be a small dg category and $n \geq 0$. Suppose that the dg functor $\ca \longrightarrow 0$ (where $0$ denotes the terminal object in $\dgcat$) has the right lifting property with respect to the set $\{ S(m)\,|\, m>n\}$. Then for all objects $x,y \in \ca$, the homology $R$-modules $\mathsf{H}_i(\ca(x,y))$ are zero for $i >n$.
\end{lemma}
\begin{proof}
This follows easily from the above definitions.
\end{proof}
\begin{lemma}\label{faltava}
Let $\xymatrix{\pi: M_{\bullet} \ar@{->>}[r]& N_{\bullet}}$ be a fibration in $Ch$ and $n+1 >0$. If the induced map on the homology $R$-modules
$$ \mathsf{H}_i(M_{\bullet}) \stackrel{\sim}{\longrightarrow} \mathsf{H}_i(N_{\bullet})$$
is an isomorphism for $i >n+1$, then $\pi$ has the right lifting property with respect to the set $\{ S_m \rightarrow D_{m+1}\,|\, m>n+1\}$, see definition~\ref{novadef}. 
\end{lemma}
\begin{proof}
Consider the short exact sequence of complexes
$$\xymatrix{ 0 \ar[r] & *+<2ex>{K_{\bullet}} \ar@{>->}[r]^i & M_{\bullet} \ar@{->>}[r]^{\pi} & N_{\bullet} \ar[r] & 0\,,}$$
where $K_{\bullet}$ denotes the kernel of $\pi$. Notice that in the induced long exact sequence on homology, the isomorphisms $(m>n+1)$
$$ \cdots \rightarrow \mathsf{H}_{m+1}(M_{\bullet}) \stackrel{\sim}{\longrightarrow} \mathsf{H}_{m+1}(N_{\bullet}) \longrightarrow \mathsf{H}_{m}(K_{\bullet}) \longrightarrow \mathsf{H}_{m}(M_{\bullet}) \stackrel{\sim}{\longrightarrow} \mathsf{H}_{m}(N_{\bullet})\longrightarrow \cdots\,,$$
imply that $\mathsf{H}_m(K_{\bullet})=0$. Now a simple diagram chasing argument (see~\cite[2.3.5]{HoveyLivro}) allow us to conclude the proof.
\end{proof}
\begin{corollary}\label{faltava1}
Let $F:\ca \rightarrow \cb$ be a dg functor such that for all objects $x,y \in \ca$, the induced morphism
$$ F(x,y): \ca(x,y) \longrightarrow \cb(Fx,Fy)$$
in $Ch$ satisfies the conditions of lemma~\ref{faltava}. Then $F$ has the right lifting property with respect to the elements of the set $\{S(m)\,|\, m>n+1\}$.
\end{corollary}

Now, let $\ca$ be a positively graded dg category. For each $n \geq 0$, apply the small object argument (\cite[10.5.14]{Hirschhorn}) to the dg functor $\ca \longrightarrow 0$, using the set $\{S(m)\,|\,m > n\}$ of generating cofibrations (see~\ref{novadef}). We obtain the following factorization
$$
\xymatrix{
\ca \ar[rr] \ar[dr]_{P_n} && 0 \\
& P_n(\ca) \ar[ur] & \,,
}
$$
where the dg functor $P_n$ is obtained by an infinite composition of pushouts along the elements of the set $\{S(m)\,|\,m > n\}$. Notice that the small object argument furnishes us natural dg functors $P_{n+1}(\ca) \longrightarrow P_n(\ca)$ making the following diagram
$$
\xymatrix{
& &\vdots \ar[d] \\
& &P_2(\ca) \ar[d] \\
& &P_1(\ca) \ar[d] \\
\ca \ar[rr]^{P_0} \ar[urr]^{P_1} \ar[uurr]^{P_2} & &P_0(\ca)
}
$$
commutative. Moreover, by construction, all the dg functors in the diagram induce the identity map on the set of objects.
\begin{proposition}\label{proofBig}
The above construction is a Postnikov tower for $\ca$.
\end{proposition}
\begin{proof}
We verify the conditions of definition~\ref{maindef}:
\begin{itemize}
\item[A1)] Since $P_n(\ca)$ is obtained by an infinite composition of pushouts along the elements of the set $\{ S(m)\,|\, m>n\}$ and the homology functors commute with infinite compositions, it is enough to prove the following: let $\cb$ be a positively graded dg category and consider the following pushout ($m>n$)
$$
\xymatrix{
\cc(m) \ar[d]_{S(m)} \ar[r]^T \ar@{}[dr]|{\lrcorner} & \cb \ar[d] \\
\cd(m+1) \ar[r] & \widetilde{\cb}
}
$$
in $\dgcat$. We need to show that $\widetilde{\cb}$ is also positively graded and that for all objects $x,y \in \cb$, the induced map on the homology $R$-modules
$$ \mathsf{H}_i(\cb(x,y)) \stackrel{\sim}{\longrightarrow} \mathsf{H}_i(\widetilde{\cb}(x,y))$$
is an isomorphism for $i \leq n$.

Observe that, as in Drinfeld's description of the Hom complexes in his DG quotient~\cite[3.1]{Drinfeld}, we have an isomorphism of graded $R$-modules (but not an isomorphism of complexes)
$$\widetilde{\cb}(x,y) \stackrel{\sim}{\longrightarrow} \bigoplus_{l=0}^{\infty} \widetilde{\cb}^l(x,y)\,,$$
where
$\widetilde{\cb}^l(x,y)$ is by definition the graded $R$-module
$$\underbrace{\cb(T(2), y)\otimes R[m+1]\otimes \cdots
\otimes  \cb(T(2), T(1)) \otimes R[m+1]
\otimes \cb(x,T(1))}_{l \textrm{ {\scriptsize factors} } R[m+1]}.
$$
The differential of an element 
$$\underbrace{g_{n+1}\cdotp h \cdots g_2 \cdotp h  \cdotp g_1}_{l \textrm{ {\scriptsize factors} } h} \in  \widetilde{\cb}^l(x,y)$$
is equal to
$$
d(g_{n+1})\cdotp h \cdots g_2 \cdotp h  \cdotp g_1 + \underbrace{(-1)^{\mid g_{n+1} \mid}\cdotp
   g_{n+1}\cdotp d(h)  \cdots g_2 \cdotp h  \cdotp  g_1}_{(l-1) \textrm{ {\scriptsize factors} } h} + \cdots \;\; \,,
$$
where $d(h) \in \cb(T(1),T(2))$ corresponds to the image of ${\bf 1}_m \in S_m$ (see~\ref{novadef}) under the dg functor $T$. This implies that, for every $j \geq 0$, the sum
$$
\bigoplus_{l \geq 0}^j \widetilde{\cb}^l(x,y) \hookrightarrow \widetilde{\cb}(x,y)
$$
is a subcomplex and so we obtain an exhaustive filtration of $\widetilde{\cb}(x,y)$. Observe that since $m>n$ and $\cb$ is positively graded the natural inclusion
$$\cb(x,y)= \widetilde{\cb}^0(x,y) \hookrightarrow \widetilde{\cb}(x,y)$$
induces isomorphisms
$$ \cb(x,y)_i \stackrel{\sim}{\longrightarrow} \widetilde{\cb}(x,y)_i$$
for $i \leq n+1$ and so an isomorphism
$$ \tau_{\leq n}\cb(x,y) \stackrel{\sim}{\longrightarrow} \tau_{\leq n}\widetilde{\cb}(x,y)$$
between the truncated complexes. We conclude that $\widetilde{\cb}$ is positively graded and that the induced map on the homology $R$-modules
$$ \mathsf{H}_i(\cb(x,y)) \stackrel{\sim}{\longrightarrow} \mathsf{H}_i(\widetilde{\cb}(x,y))$$
is an isomorphism for $i \leq n$.
\item[A2)] By condition $A1)$, for all objects $x,y \in \ca$, the induced map on the homology $R$-modules
$$ \mathsf{H}_0(\ca(x,y)) \stackrel{\sim}{\longrightarrow} \mathsf{H}_0(P_n(\ca)(x,y))$$
is an isomorphism. Since the dg functor $P_n: \ca \longrightarrow \ca_n$ induces the identity map on the set of objects, we conclude that the induced functor $\mathsf{H}_0(\ca) \stackrel{\sim}{\longrightarrow}\mathsf{H}_0(\ca_n)$ is an equivalence of categories.
\item[B)] By construction, the dg functor $P_n(\ca) \longrightarrow 0$ has the right lifting property with respect to the set $\{ S(m)\,|\,m>n\}$. This implies, by lemma~\ref{lift}, that for all objects $x,y \in \ca$ the homology $R$-modules $\mathsf{H}_i(P_n(\ca)(x,y))$ are zero for $i>n$.
\end{itemize}
\end{proof}
\begin{notation}
We denote by $P(\ca)$ the Big Postnikov model thus obtained.
\end{notation}

\subsection{Uniqueness and homotopy type}

\begin{proposition}\label{existence}
Let $(\ca_n)_{n \geq 0}$ be a Postnikov tower for a homologically connective dg category $\ca$, where all the dg functors $\xymatrix{\ca_{n+1} \ar@{->>}[r] & \ca_n}$ are fibrations. Then there exists a morphism
$$ M: P(\ca) \longrightarrow (\ca_n)_{n\geq 0}$$
between Postnikov towers.
\end{proposition}
\begin{proof}
We will construct $M$ recursively, starting with the case $n=0$ and then going upwards.

($n=0$) Notice that the small object argument allows us to construct inductively a dg functor $M_0: P_0(\ca) \longrightarrow \ca_0$ as follows:

\vspace{0.2cm}

{\it step:} suppose we have the following (solid) diagram $(i \geq 0, \,P_0(\ca)^0=\ca)$
$$
\xymatrix{
\underset{m>0}{\coprod} \,\, \underset{\cc(m) \rightarrow P_0(\ca)^i}{\coprod} \cc(m) \ar[rr] \ar@{}[drr]|{\lrcorner} \ar[d] \ar@/^2pc/[rrrr]^{T_i} && P_0(\ca)^i \ar[d] \ar[rr]_{M^i_0} & &\ca_0 \\
\underset{m>0}{\coprod} \,\, \underset{\cc(m) \rightarrow P_0(\ca)^i}{\coprod} \cd(m+1) \ar[rr] & &P_0(\ca)^{i+1} \ar@{-->}[urr]_{M_0^{i+1}}& &\,.
}
$$
Recall from that we denote by ${\bf 1}_m$ the cycle of degree $m$ in $S_m$ (and so in $\cc(m)(1,2)$) which corresponds to the unit of $R$. Since $\ca_0$ satisfies condition $B)$ of definition~\ref{maindef}, we can choose a bounding chain $b$ in $\ca_0$ for each cycle $T_i({\bf 1}_m), \,m>0$ (i.e. $d(b)=T_i({\bf 1}_m)$). These choices give rise to a dg functor $M_0^{i+1}$ which makes the above diagram commute.

\vspace{0.2cm}

By passing to the colimit on $i$, we obtain our desired dg functor
$$M_0=\underset{i}{\mbox{colim}}\,M_0^i: P_0(\ca)= \underset{i}{\mbox{colim}}\,P_0(\ca)^i \longrightarrow \ca_0\,.$$

$(n \Rightarrow n+1)$ Suppose we have a dg functor $M_n:P_n(\ca) \longrightarrow \ca_n$ between the $n$th Postnikov sections. We will construct a `lift' $M_{n+1}$ which makes the square
$$
\xymatrix{
P_{n+1}(\ca) \ar[d] \ar@{-->}[rr]^{M_{n+1}} & &\ca_{n+1} \ar@{->>}[d]\\
P_n(\ca) \ar[rr]_{M_n} & &\ca_n
}
$$
commutative. Our argument is also an inductive one:

\vspace{0.2cm}

{\it step:} suppose we have the following (solid) diagram $(i \geq 0,\, P_{n+1}(\ca)^0=\ca)$
$$
\xymatrix{
\underset{m > n+1}{\coprod} \,\, \underset{\cc(m) \rightarrow P_{n+1}(\ca)^i}{\coprod} \cc(m) \ar[rr] \ar[d] \ar@{}[drr]|{\lrcorner} & & P_{n+1}(\ca)^i \ar[d] \ar[rr]^{M^i_{n+1}} && \ca_{n+1} \ar@{->>}[d] \\
\underset{m > n+1}{\coprod} \,\, \underset{\cc(m) \rightarrow P_{n+1}(\ca)^i}{\coprod} \cd(m+1) \ar[rr] && P_{n+1}(\ca)^{i+1} \ar@{-->}[rr]_{\widetilde{M^{i+1}_n}} \ar@{-->}[urr]_{M^{i+1}_{n+1}} & &\ca_n \,.
}
$$
Notice that the left (solid) square appears in the construction of $P_n(\ca)^{i+1}$. This implies that the dg functor $M_n^{i+1}: P_n(\ca)^{i+1} \longrightarrow \ca_n$ restricts to a dg functor $\widetilde{M_n^{i+1}}$, which makes the right square commutative. Now, observe that the dg functor $\xymatrix{\ca_{n+1} \ar@{->>}[r] &\ca_n}$ satisfies the conditions of corollary~\ref{faltava1} and so it has the right lifting property with respect to the elements of the set $\{ S(m)\,|\, m>n\}$. This implies that there exists an induced dg functor $M_{n+1}^{i+1}$ which makes the above diagram commute.

\vspace{0.2cm}

By passing to the colimit on $i$, we obtain our desired morphism
$$M_{n+1}=\underset{i}{\mbox{colim}}\,M_{n+1}^i: P_{n+1}(\ca)= \underset{i}{\mbox{colim}}\,P_{n+1}(\ca)^i \longrightarrow \ca_{n+1}\,.$$

The proof is now finished.
\end{proof}
\begin{remark}\label{big->small}
Since in the small Postnikov model $\mathbb{P}(\ca)$ for $\ca$, the dg functors
$$\xymatrix{ \tau_{\leq n+1}(\ca) \ar@{->>}[r] & \tau_{\leq n}(\ca)}$$
are fibrations, proposition~\ref{existence} implies the existence of a morphism
$$ M: P(\ca) \longrightarrow \mathbb{P}(\ca)$$
from the Big to the small Postnikov model. Moreover, the bounding chains in $\mathbb{P}(\ca)$ used in the construction of $M$ are all trivial and so this morphism is well-defined. 
Notice also that for $n \geq 0$, the dg functor $M_n$ satisfies all the conditions of corollary~\ref{corfib} and so it is a fibration in $\dgcat$. 
\end{remark}
We now prove that Postnikov towers are `essentially' unique.
\begin{theorem}\label{uniqueness}
Let $\ca$ be a homologically connective dg category. Given two objects in $\mbox{Post}(\ca)$ (see~\ref{catPost}), there exists a zig-zag of weak equivalences (see~\ref{morph}) relating the two.
\end{theorem}
\begin{proof}
Let $(\ca_n)_{n \geq 0}$ and $(\ca_n')_{n \geq 0}$ two Postnikov towers for $\ca$. By remark~\ref{fibreplace}, we can construct morphisms in $\mbox{Post}(\ca)$
$$ \xymatrix{(\ca_n)_{n \geq 0} \ar[r]^-{\sim} & F(\ca_n)_{n \geq 0} & (\ca_n')_{n \geq 0} \ar[r]^-{\sim} & F(\ca_n')_{n \geq 0}\,,}$$
such that the dg functors
$$ \xymatrix{F(\ca_{n+1}) \ar@{->>}[r] & F(\ca_n) & F(\ca_{n+1}') \ar@{->>}[r] & F(\ca_n')}$$
are fibrations in $\dgcat$. Moreover, by proposition~\ref{existence}, we can also construct morphisms as follows
$$ \xymatrix{P(\ca) \ar[r]^-{\sim} & F(\ca_n)_{\geq 0} & P(\ca) \ar[r]^-{\sim} & F(\ca_n')_{n \geq 0}\,.}$$
We obtain finally, the following zig-zag
$$ \xymatrix{ (\ca_n)_{n \geq 0} \ar[r]^-{\sim} & F(\ca_n)_{n \geq 0}  & P(\ca) \ar[l]_-{\sim} \ar[r]^-{\sim} & F(\ca_n')_{n \geq 0} & (\ca_n')_{n \geq 0}\ar[l]_-{\sim}}$$
of weak equivalences in $\mbox{Post}(\ca)$. 
\end{proof}
\begin{remark}
Notice that by theorem~\ref{uniqueness}, the classifying space (\cite[14]{Hirschhorn}) of $\mbox{Post}(\ca)$ has a single connected component.
\end{remark}

We now show how the full homotopy type of a homologically connective dg category can be recovered from any of its Postnikov towers.
\begin{proposition}\label{homtype}
Let $\ca$ be a homologically connective dg category and $(\ca_n)_{n \geq 0}$ a Postnikov tower for $\ca$. Then the natural dg functor
$$ \ca \longrightarrow \underset{n}{\mbox{holim}}\, \ca_n$$
is a quasi-equivalence.
\end{proposition}
\begin{proof}
Notice that theorem~\ref{uniqueness} and remark~\ref{morph} imply that the homotopy limit of any Postnikov tower for $\ca$ is well defined up to quasi-equivalence. We can then consider the small Postnikov model $\mathbb{P}(\ca)$ for $\ca$. Since every object in $\dgcat$ is fibrant (see~\ref{objfib}) and the dg functors
$$\xymatrix{\tau_{\leq n+1}(\ca) \ar@{->>}[r] & \tau_{\leq n}(\ca)}$$
are fibrations in $\dgcat$, we have a natural quasi-equivalence
$$ \underset{n}{\mbox{lim}}\, \tau_{\leq n}(\ca) \stackrel{\sim}{\longrightarrow} \underset{n}{\mbox{holim}}\, \tau_{\leq n}(\ca)\,.$$
By construction of limits in $\dgcat$, we conclude that the natural dg functor
$$\ca \stackrel{\sim}{\longrightarrow} \underset{n}{\mbox{lim}}\,\mathbb{P}_n(\ca)$$
is an isomorphism.

\end{proof}
\section{$k$-invariants}
In this chapter we construct $k$-invariants for homologically connective dg categories (see definitions~\ref{def1} and \ref{def2}). We show that these invariants correspond to derived derivations with values in a certain bimodule (see~\ref{Dderivations}). Then we prove our main theorem (\ref{mainthm}), which shows how the full homotopy type of the $n+1$ Postnikov section of an homologically connective dg category $\ca$ can be recovered from the $n$th $k$-invariant of $\ca$. For constructions of $k$-invariants in the context of spectral algebra see~\cite{Dugger}~\cite{DS}~\cite{Lazarev}.
Let us start with some general constructions.

\begin{definition}\label{squarezero}
Let $\ca$ be a small dg category and $M$ a $\ca\text{-}\ca$-bimodule (see~\ref{bimodule}). The {\em square zero extension $\ca \ltimes M$ of $\ca$ by $M$} is the dg category defined as follows: its objects are those of $\ca$ and for objects $x,y \in \ca$ we have
$$ \ca \ltimes M(x,y):= \ca(x,y) \oplus M(x,y)\,.$$
The composition in $\ca \ltimes M$ is defined using the composition on $\ca$, the above bimodule structure and by imposing that the composition between $M$-factors is zero. 
\end{definition}
\begin{remark}\label{rksquare0}
Notice that $\ca$ is a (non-full) dg subcategory of $\ca \ltimes M$ and that we have a natural projection dg functor
$$ \xymatrix{\ca \ltimes M \ar@{->>}[r] & \ca}\,,$$
which is clearly a fibration in $\dgcat$, see proposition~\ref{deffibrations}.
\end{remark}
\begin{definition}\label{derivations}
\begin{itemize}
\item[-] A {\em derivation of $\ca$ with values in a $\ca\text{-}\ca$-bimodule $M$} is a morphism in $\dgcat\downarrow \ca$ (see~\ref{overcat}) from $\ca$ to $\ca \ltimes M$, or equivalently a section of the natural projection dg functor $\xymatrix{ \ca \ltimes M \ar@{->>}[r] & \ca}$.
\item[-] A {\em derived derivation of $\ca$ with values in a $\ca\text{-}\ca$-bimodule $M$} is a morphism in the homotopy category $\mathsf{Ho}(\dgcat\downarrow \ca)$ (see~\ref{rks}) from $\ca$ to $\ca \ltimes M$.
\end{itemize}
\end{definition}
\begin{notation}\label{hoje1}
We denote by $\mbox{Der}(\ca, M)$ (resp. $\mathbb{R}\mbox{Der}(\ca,M)$) the set of derivations (resp. derived derivations) of $\ca$ with values in $M$. The (derived) derivation obtained by considering $\ca$ as a dg subcategory of $\ca \ltimes M$ is called the trivial one.
\end{notation}
\begin{remark}
Notice that if $\ca$ is a $R$-algebra $A$ (i.e. $\ca$ has only one object and its endomorphisms $R$-algebra is $A$), the notion of derivation coincides with the classical one, i.e. a $R$-linear map $D:A \longrightarrow M$ which satisfies the  Leibniz relation
$$ D(ab)=a(Db)+(Da)b\,\,\,\,\,\,\,\, a,b \in A\,.$$
\end{remark}
\begin{proposition}\label{hoje}
Let $F:\ca \longrightarrow \cb$ be an object in $\dgcat \downarrow \cb$ and $M$ a $\cb\text{-}\cb$-bimodule. Then the set $\mathsf{Ho}(\dgcat \downarrow \cb)(\ca, \cb \ltimes M)$ is naturally isomorphic to the set of derived derivations $\mathbb{R}\mbox{Der}(\ca, F^{\ast}(M))$ of $\ca$ with values in the $\ca\text{-}\ca$-bimodule $F^{\ast}(M)$ obtained by restricting $M$ along $F$.
\end{proposition}
\begin{proof}
Recall from remark~\ref{rks}, the (derived) Quillen adjunction
$$
\xymatrix{
\mathsf{Ho}(\dgcat \downarrow \ca) \ar@<-1ex>[d]_{F_!} \\
\mathsf{Ho}(\dgcat \downarrow \cb) \ar@<-1ex>[u]_{\mathbb{R}F^!}\,.
}
$$
Notice that we have the following pull-back square
$$
\xymatrix{
\ca \ltimes F^{\ast}(M) \ar@{->>}[d] \ar[rr]^{F \ltimes Id} \ar@{}[drr]|{\ulcorner} && \cb \ltimes M \ar@{->>}[d] \\
\ca \ar[rr]_F && \cb\,,
}
$$
which shows us that the image of $\cb \ltimes M$ under the functor $\mathbb{R}F^!$ is isomorphic to $\ca \ltimes F^{\ast}(M)$. Moreover the image of $\ca$ under the functor $F_!$ is isomorphic to the object $F:\ca \longrightarrow \cb$ in $\mathsf{Ho}(\dgcat\downarrow \cb)$ and so by adjunction we obtain the desired isomorphism.
\end{proof}

We now define the dg categories which play the same role as the Eilenberg-Mac Lane spaces in the classical theory of $k$-invariants.
\begin{definition}\label{Eilenberg}
Let $\ca$ be a positively graded dg category and $n \geq 0$. Consider the following bimodule:
$$ 
\begin{array}{rcl}
\mathsf{H}_{n+1}(\ca)[n+2]: \mathsf{H}_0(\ca)^{op}\otimes \mathsf{H}_0(\ca) & \longrightarrow & Ch \\
(x,y) & \mapsto & \mathsf{H}_{n+1}((\ca)(x,y))[n+2]\,,
\end{array}
$$
where the complex $\mathsf{H}_{n+1}((\ca)(x,y))[n+2]$ is simply the $R$-module $\mathsf{H}_{n+1}(\ca(x,y))$ concentrated in degree $n+2$. Notice that the natural projection dg functor $\mathbb{P}_n(\ca) \longrightarrow \mathbb{P}_0(\ca)=\mathsf{H}_0(\ca)$ endow $\mathsf{H}_{n+1}(\ca)[n+2]$ with a structure of $\mathbb{P}_n(\ca)\text{-}\mathbb{P}_n(\ca)$-bimodule. Finally, we denote by $\ca \ltimes \mathsf{H}_{n+1}(\ca)[n+2]$ the square zero extension obtained (\ref{squarezero}) using this bimodule structure.
\end{definition}
\begin{remark}\label{natpoint}
Notice that by remark~\ref{rksquare0}, $\mathbb{P}_n(\ca)$ is a dg subcategory of $\mathbb{P}_n(\ca) \ltimes \mathsf{H}_{n+1}(\ca)[n+2]$ and we have a natural projection dg functor
$$ \xymatrix{ \mathbb{P}_n(\ca) \ltimes~\mathsf{H}_{n+1}(\ca)[n+~2] \ar@{->>}[r]& \mathbb{P}_n(\ca)\,.}$$
\end{remark}
\begin{definition}\label{k-invariant}
Let 
$$ \gamma_n: P_n(\ca) \longrightarrow \mathbb{P}_n(\ca) \ltimes \mathsf{H}_{n+1}(\ca)[n+2]$$ be the natural dg functor obtained by modifying the dg functor $\xymatrix{M_n: P_n(\ca) \ar@{->>}[r] & \mathbb{P}_n(\ca)}$ (see \ref{big->small}) as follows:

\vspace{0.2cm}

{\it step:} suppose we have the following (solid) diagram $(i \geq 0, \,P_n(\ca)^0=\ca)$
$$
\xymatrix{
\underset{m > n}{\coprod} \,\, \underset{\cc(m) \rightarrow P_n(\ca)^i}{\coprod} \cc(m) \ar[rr]^-{T_i} \ar[d] \ar@{}[drr]|{\lrcorner} && P_n(\ca)^i \ar[rr]^-{\gamma_n^i} \ar[d]  && \mathbb{P}_n(\ca)\ltimes \mathsf{H}_{n+1}(\ca)[n+2] \\
\underset{m > n}{\coprod} \,\, \underset{\cc(m) \rightarrow P_n(\ca)^i}{\coprod} \cd(m+1) \ar[rr] && P_{n}(\ca)^{i+1} \ar@{-->}[urr]_{\gamma_n^{i+1}} & & \,.
}
$$
For every cycle $T_i(\mathbf{1}_m), \, m >n+1$ choose $0$ as a bounding chain in $\mathbb{P}_n(\ca)$ (and so in $\mathbb{P}_n(\ca) \ltimes \mathsf{H}_{n+1}(\ca)[n+2]$), as in the case of the dg functor $M_n$. Now, let $T_i(\mathbf{1}_{n+1}) \in P_n(\ca)^i(T_i(1), T_i(2))$ be a cycle of degree $n+1$. Since $m>n$, the description of the complexes of morphisms in $P_n(\ca)^i$ (see proof of proposition~\ref{proofBig}) implies that we have natural isomorphisms
$$ \ca(T_i(1),T_i(2))_j \stackrel{\sim}{\longrightarrow} P_n(\ca)^i(T_i(1), T_i(2))_j$$
for $j \leq n+1$. We can then choose for bounding chain for $T_i({\bf 1}_{n+1})$ its homology class in $\mathsf{H}_{n+1}(\ca(T_i(1), T_i(2))$. These choices give rise to a dg functor $\gamma_n^{i+1}$ which makes the above diagram commute.

\vspace{0.2cm}

By passing to the colimit on $i$, we obtain our desired dg functor
$$ \gamma_n= \underset{i}{\mbox{colim}}\, \gamma_n^i: P_n(\ca)= \underset{i}{\mbox{colim}}\, P_n(\ca)^i \longrightarrow \mathbb{P}_n(\ca)\ltimes \mathsf{H}_{n+1}(\ca)[n+2]\,.$$
\end{definition}
\begin{remark}\label{gammafib}
Notice that by construction, the dg functor $\gamma_n$ satisfies all the conditions of corollary~\ref{corfib} and so it is a fibration in $\dgcat$. Moreover for $n \geq 0$, we have the following commutative diagram in $\dgcat$
$$
\xymatrix{
P_n(\ca) \ar@{->>}[dr]^{\sim}_{M_n} \ar@{->>}[rr]^-{\gamma_n} && \mathbb{P}_n(\ca) \ltimes \mathsf{H}_{n+1}(\ca)[n+2] \ar@{->>}[dl] \\
& \mathbb{P}_n(\ca) & \,.
}
$$
\end{remark}
\begin{notation}\label{important}
We denote by $\dgcat \downarrow \mathbb{P}_n(\ca)$ the category of objects in $\dgcat$ over $\mathbb{P}_n(\ca)$, see notation~\ref{overcat}. 
\end{notation}
\begin{definition}\label{def1}
Let $\ca$ be a positively graded dg category and $n \geq 0$. Its {\em $n$th $k$-invariant $\alpha_n(\ca)$} is by definition the image of the dg functor $\gamma_n$ in the homotopy category $\mathsf{Ho}(\dgcat\downarrow \mathbb{P}_n(\ca))$, see remark~\ref{gammafib}
\end{definition}
\begin{remark}\label{Dderivations}
Since the dg functor $\xymatrix{M_n: P_n(\ca) \ar@{->>}[r]^-{\sim} & \mathbb{P}_n(\ca)}$ is a quasi-equivalence, we have an isomorphism between
$$ \mathsf{Ho}(\dgcat\downarrow \mathbb{P}_n(\ca))(P_n(\ca), \mathbb{P}_n(\ca) \ltimes \mathsf{H}_{n+1}(\ca)[n+2])$$ and $$\mathsf{Ho}(\dgcat\downarrow \mathbb{P}_n(\ca))(\mathbb{P}_n(\ca), \mathbb{P}_n(\ca) \ltimes \mathsf{H}_{n+1}(\ca)[n+2])$$
which implies that $\alpha_n(\ca)$ corresponds to a derived derivation of $\mathbb{P}_n(\ca)$ with values in the $\mathbb{P}_n(\ca)\text{-}\mathbb{P}_n(\ca)$-bimodule $\mathsf{H}_{n+1}(\ca)[n+2]$, see definiton~\ref{derivations}.
\end{remark}
\begin{definition}\label{def2}
Let $\ca$ be a homologically connective dg category. Its $n$th $k$-invariant $\alpha_n(\ca)$ is by definition the $n$th $k$-invariant of $\tau_{\geq 0}(\ca)$, see remark~\ref{pos->con}.
\end{definition}
\begin{remark}\label{natpoint1}
Notice that although the category $\dgcat \downarrow \mathbb{P}_n(\ca)$ is not pointed (the initial and terminal objects are not isomorphic), there is a natural morphism (in $\dgcat \downarrow \mathbb{P}_n(\ca)$) from its terminal object $\mathbb{P}_n(\ca)$ to $\mathbb{P}_n(\ca) \ltimes \mathsf{H}_{n+1}(\ca)[n+2]$ (see remark~\ref{natpoint}).
\end{remark}
We now show how the full homotopy type of $P_{n+1}(\ca)$ in $\dgcat$ can be entirely recovered from the $n$th $k$-invariant $\alpha_n(\ca)$.
\begin{theorem}\label{mainthm}
We have a homotopy fiber sequence
$$\xymatrix{ P_{n+1}(\ca) \ar[r] & P_n(\ca) \ar@{->>}[r]^-{\gamma_n} & \mathbb{P}_n(\ca)\ltimes \mathsf{H}_{n+1}(\ca)[n+2]}$$
in $\mathsf{Ho}(\dgcat\downarrow \mathbb{P}_n(\ca))$.
\end{theorem}
\begin{proof}
We need to show that $P_{n+1}(\ca)$ is quasi-equivalent in $\dgcat$ to the homotopy pullback of the diagram
$$
\xymatrix{
&& \mathbb{P}_n(\ca) \ar[d] \\
P_n(\ca) \ar@{->>}[rr]_-{\gamma_n} && \mathbb{P}_n(\ca) \ltimes \mathsf{H}_{n+1}(\ca)[n+2]\,.
}
$$
Since $\gamma_n$ is a fibration (see~\ref{gammafib}) and every dg category is fibrant (see~\ref{objfib}), the homotopy pullback and the pullback are quasi-equivalent. Notice that we have the following commutative diagram
$$
\xymatrix@!0 @R=3.5pc @C=7pc{
& \ca \ar[d]_{P_{n+1}} \ar@/_2pc/[dddl]_{P_n} \ar@/^2pc/[dddr] &  \\
& P_{n+1}(\ca) \ar@/_1pc/[ddl] \ar@/^1pc/[ddr] &  \\
& \cw \ar[dr] \ar[dl] & \\
P_n(\ca) \ar@{->>}[dr]_-{\gamma_n} && \mathbb{P}_n(\ca)\ar[dl] \\
& \mathbb{P}_n(\ca) \ltimes \mathsf{H}_{n+1}(\ca)[n+2] \ar@{}[uu]|-{\widehat{}} &\,.
}
$$
This diagram gives rise to the following factorization
$$
\xymatrix{
\ca \ar@{-->}[drr]_{\phi} \ar[rr]^-{P_{n+1}} && P_{n+1}(\ca) \ar@{-->}[d]^-{\theta}\\
&&  \cw\,,
}
$$
where $\theta$ and $\phi$ are the induced dg functors to the pullback $\cw$. We need to show that $\theta$ is a quasi-equivalence. By construction of limits in $\dgcat$, all the dg functors in the previous diagrams induce the identity map on the set of objects and so it is enough to prove that for all objects $x,y \in P_{n+1}(\ca)$, the morphism of complexes
$$ \theta(x,y): P_{n+1}(\ca)(x,y) \longrightarrow \cw(x,y)$$
is a quasi-isomorphism. Let us denote by
$$ 0 \leftarrow M_0 \leftarrow M_1 \leftarrow \cdots \leftarrow M_n \leftarrow M_{n+1} \leftarrow M_{n+2} \leftarrow M_{n+3} \leftarrow \cdots$$
the complex $\ca(x,y)$. Notice that by construction of $P_n(\ca)$ (see~\ref{proofBig}), the complex $P_n(\ca)(x,y)$ is of the following shape
$$ 0 \leftarrow M_0 \leftarrow M_1 \leftarrow \cdots \leftarrow M_n \leftarrow M_{n+1} \leftarrow \widetilde{M_{n+1}} \leftarrow \widetilde{M_{n+3}} \leftarrow \cdots\,.$$
The complex $\cw(x,y)$ identifies then with the pullback of the following diagram
$$
\xymatrix@!0 @R=2.5pc @C=8pc{
\vdots \ar[d] & \vdots \ar[d]  & \vdots \ar[d] \\
\widetilde{M_{n+3}} \ar[r] \ar[d] & 0 \ar[d] \ar@{=}[r] & 0 \ar[d]\\
\widetilde{M_{n+2}} \ar@{->>}[r] \ar[d] & \mathsf{H}_{n+1}(\ca(x,y)) \ar[d]  & 0 \ar[d] \ar[l]\\
M_{n+1} \ar[r] \ar[d] & 0 \ar[d] \ar@{=}[r] & 0 \ar[d]\\
M_n \ar@{->>}[r] \ar[d] & M_n/\mbox{Im}(M_{n+1}) \ar[d] \ar@{=}[r] & M_n/\mbox{Im}(M_{n+1}) \ar[d]\\
\vdots \ar[d] & \vdots \ar[d] & \vdots \ar[d]\\
M_1 \ar@{=}[r] \ar[d] & M_1 \ar@{=}[r] \ar[d] & M_1 \ar[d]\\
M_0 \ar@{=}[r] \ar[d] & M_0 \ar@{=}[r] \ar[d] & M_0 \ar[d]\\
0 & 0 &0 \,.
}
$$
The above diagram allow us to conclude that $\mathsf{H}_j \cw(x,y)=0$ for $j \geq n+2$ and that the induced map
$$ \mathsf{H}_j(P_{n+1}(\ca)(x,y)) \stackrel{\sim}{\longrightarrow} \mathsf{H}_j\cw(x,y)$$
is an isomorphism for $j \neq n+1$.

We now prove that the induced map
$$ \mathsf{H}_{n+1}\ca(x,y) \stackrel{\sim}{\longrightarrow} \mathsf{H}_{n+1}\cw(x,y)$$
is an isomorphism. Notice that this implies (by the 2 out of 3 property) that $\theta(x,y)$ is a quasi-isomorphism. In order to prove this, we start by observing that in $\dgcat$, pullbacks commute with filtered colimits. Since $P_n(\ca)$ is constructed as a filtered colimit and the homology functor $\mathsf{H}_{n+1}(-)$ preserves filtered colimits it is then enough to prove the following: 

\vspace{0.2cm}

{\it start:} consider the following pullback square
$$
\xymatrix{
\ca \ar@{=}[d] \ar[rr] \ar@{}[drr]|{\ulcorner} && \mathbb{P}_n(\ca) \ar[d] \\
\ca \ar[rr] && \mathbb{P}_n(\ca) \ltimes \mathsf{H}_{n+1}(\ca)[n+2]\,.
}
$$

\vspace{0.2cm}

{\it step:} consider the commutative diagram ($i \geq 0, \,P_n(\ca)^0=\ca$)
$$
\xymatrix{
 \cc(n+1) \ar[r]^-{T} \ar[d]_{S(n+1)} \ar@{}[dr]|{\lrcorner} & P_n(\ca)^i \ar[rr] \ar[d] & & \mathbb{P}_n(\ca)\ltimes \mathsf{H}_{n+1}(\ca)[n+2] \\
 \cd(n+2) \ar[r] & \widetilde{P_{n}(\ca)^{i}} \ar@{-->}[urr]_{\gamma_{i}(T)} &
}
$$
used in the construction of the natural dg functor $\gamma_n$ (see \ref{k-invariant}), and suppose that the induced dg functor from $\ca$ to the pullback
$$
\xymatrix{
\cw_{i}(T) \ar[d] \ar[rr] \ar@{}[drr]|{\ulcorner} && \mathbb{P}_n(\ca) \ar[d] \\
P_n(\ca)^i \ar[rr] && \mathbb{P}_n(\ca) \ltimes \mathsf{H}_{n+1}(\ca)[n+2]
}
$$
induces an isomorphism
$$ \mathsf{H}_{n+1}\ca(x,y) \stackrel{\sim}{\longrightarrow} \mathsf{H}_{n+1}(\cw_{i}(T)(x,y))\,.$$
We need to show that the induced dg functor from $\cw_{i}(T)$ to the pullback
$$
\xymatrix{
\widetilde{\cw_{i}(T)} \ar[d] \ar[rr] \ar@{}[drr]|{\ulcorner} && \mathbb{P}_n(\ca) \ar[d] \\
\widetilde{P_n(\ca)^i} \ar[rr]_-{\gamma_i(T)} && \mathbb{P}_n(\ca) \ltimes \mathsf{H}_{n+1}(\ca)[n+2]
}
$$
induces an isomorphism
$$ \mathsf{H}_{n+1}(\cw_{i}(T)(x,y)) \stackrel{\sim}{\longrightarrow} \mathsf{H}_{n+1}(\widetilde{\cw_{i}(T)}(x,y))\,.$$

Recall that for all objects $x,y \in P_n(\ca)^i$, we have an isomorphism of graded $R$-modules
$$\widetilde{P_n(\ca)^i}(x,y) \stackrel{\sim}{\longrightarrow} \bigoplus_{l=0}^{\infty} \widetilde{P_n(\ca)^i}^l(x,y)\,,$$
where $\widetilde{P_n(\ca)^i}^l(x,y)$ is the graded $R$-module
$$ \underbrace{P_n(\ca)^i(T(2), y)\otimes R[n+2]\otimes \cdots
\otimes  P_n(\ca)^i(T(2), T(1)) \otimes R[n+2]
\otimes P_n(\ca)^i(x,T(1))}_{l \textrm{ {\scriptsize factors} } R[n+2]}.
$$
The differential of an element 
$$\underbrace{g_{n+1}\cdotp h \cdots g_2 \cdotp h  \cdotp g_1}_{l \textrm{ {\scriptsize factors} } h} \in  \widetilde{P_n(\ca)^i}^l(x,y)$$
is equal to
$$
d(g_{n+1})\cdotp h \cdots g_2 \cdotp h  \cdotp g_1 + \underbrace{(-1)^{\mid g_{n+1} \mid}\cdotp
   g_{n+1}\cdotp d(h) \cdots g_2 \cdotp h  \cdotp  g_1}_{(l-1) \textrm{ {\scriptsize factors} } h} + \cdots \;\; \,,
$$
where $d(h) \in P_n(\ca)^i(T(1),T(2))$ corresponds to the image of ${\bf 1}_{n+1} \in S_{n+1}$ (see~\ref{novadef}) under the dg functor $T$. 
This description show us that the unique elements in $\widetilde{\cw_{i}(T)}(x,y)$, which eventually `destroy' the $(n+1)$-homology of the complex $\cw_{i}(T)(x,y)$ belong to the graded $R$-module
$$ \widetilde{P_n(\ca)^i}^1(x,y) = P_n(\ca)^i(T(2),y)\otimes R[n+2]\otimes P_n(\ca)^i(x,T(1))\,.$$
We now show that if $g_2 \cdotp h \cdotp g_1$ is an (homogeneous) element of degree $n+2$ in $\widetilde{P_n(\ca)^i}^1(x,y)$, whose differential 
$$g_2 \cdotp d(h) \cdotp g_1 \in (P_n(\ca)^i(x,y))_{n+1} \simeq (\cw_{i}(T)(x,y))_{n+1}$$
is non-trivial in the homology $R$-module $\mathsf{H}_{n+1}(\cw_{i}(T)(x,y))$, then the element $g_2 \cdotp h \cdotp g_1$ does not belong to $\widetilde{\cw_{i}(T)}(x,y)$. By hypothesis we have an induced isomorphism
$$ \mathsf{H}_{n+1}\ca(x,y) \stackrel{\sim}{\longrightarrow} \mathsf{H}_{n+1}(\cw_{i}(T)(x,y))$$
and so by definition~\ref{Eilenberg}, the image of $g_2 \cdotp h \cdotp g_1$ under the dg functor $\gamma_{i}(T)$ corresponds precisely to this non-trivial element in the homology $R$-module $\mathsf{H}_{n+1}\ca(x,y)$. This implies that $g_2 \cdotp h \cdotp g_1$ does not belong to the pullback complex $\widetilde{\cw_{i}(T)}(x,y)$ and so we conclude that we have an induced isomorphism
$$ \mathsf{H}_{n+1}(\cw_{i}(T)(x,y)) \stackrel{\sim}{\longrightarrow} \mathsf{H}_{n+1}(\widetilde{\cw_{i}(T)}(x,y))\,.$$

Finally, by and infinite composition procedure, we obtain the pullback $\cw$. Since the homology functor $H_{n+1}(-)$ commutes with filtered colimits, the induced map
$$ \mathsf{H}_{n+1}\ca(x,y) \stackrel{\sim}{\longrightarrow} \mathsf{H}_{n+1}\cw(x,y)$$
is an isomorphism and so we conclude that
$$ \theta(x,y): P_{n+1}(\ca)(x,y) \stackrel{\sim}{\longrightarrow} \cw(x,y)$$
is a quasi-isomorphism.
This proves the theorem.
\end{proof}
\section{Obstruction theory}
In this chapter, we develop an obstruction theory for dg categories. Our motivation comes from the examples appearing in non-abelian Hodge theory (see~\ref{example1}). We formulate the following general `rigidification' problem.

\vspace{0.2cm}

{\bf The `rigidification' problem:} Let $\ca$ be a positively graded dg category and $F_0: \cb \longrightarrow \mathsf{H}_0(\ca)$ a dg functor with values in its homotopy category, with $\cb$ a cofibrant dg category. Is there a lift $F: \cb \longrightarrow \ca$ making the diagram
$$
\xymatrix{
& \ca \ar[d]^{\tau_{\leq 0}} \\
\cb \ar[r]_-{F_0} \ar[ur]^F & \mathsf{H}_0(\ca)
}
$$
commute?

\vspace{0.2cm}

Intuitively the dg functor $F_0$ represents the `up-to-homotopy' information that one would like to rigidify, i.e. lift to the dg category $\ca$. 
\begin{remark}
Notice that if $\ca$ is a homologically connective dg category, we have a zig-zag of dg functors
$$\xymatrix{ \ca & \ar[l]_{\sim} \tau_{\geq 0}(\ca) \ar[r]^{\tau_{\leq 0}} & \mathsf{H}_0(\ca)\,.}$$
In this situation we search for a lift $\cb \longrightarrow \ca$ which factors through $\tau_{\geq 0}(\ca)$.
\end{remark}
In order to solve this problem we consider the following notion: let $\ca$ be a positively graded dg category and recall from section~\ref{secBig} its Big Postnikov model
$$
\xymatrix{
& &\vdots \ar[d] \\
& &P_2(\ca) \ar[d] \\
& &P_1(\ca) \ar[d] \\
\ca \ar[rr]^{P_0} \ar[urr]^{P_1} \ar[uurr]^{P_2} & &P_0(\ca)\,.
}
$$
\begin{definition}\label{obstruction}
Let $F:\cb \longrightarrow P_n(\ca)$ be a dg functor. Its {\em obstruction class $\omega_n(F)$} is the image of the composed dg functor (see~\ref{k-invariant})
$$ \cb \stackrel{F}{\longrightarrow} P_n(\ca) \stackrel{\gamma_n}{\longrightarrow} \mathbb{P}_n(\ca) \ltimes \mathsf{H}_{n+1}(\ca)[n+2]$$
in the homotopy category $\mathsf{Ho}(\dgcat \downarrow \mathbb{P}_n(\ca))$, see remark~\ref{gammafib}.

We say that the obstruction class $\omega_n(F)$ {\em vanishes} if it factors through the canonical morphism
$$ \mathbb{P}_n(\ca) \longrightarrow \mathbb{P}_n(\ca) \ltimes \mathsf{H}_{n+1}(\ca)[n+2]$$
in $\dgcat \downarrow \mathbb{P}_n(\ca)$, see remark~\ref{natpoint1}.
\end{definition}
\begin{remark}\label{hoje3}
Consider the composed dg functor $\xymatrix{\cb \ar[r]^-F & P_n(\ca) \ar@{->>}[r]^-{M_n}_-{\sim} & \mathbb{P}_n(\ca)}$ as an object in $\dgcat\downarrow \mathbb{P}_n(\ca)$. By proposition~\ref{hoje}, the set $$\mathsf{Ho}(\dgcat\downarrow \mathbb{P}_n(\ca))(\cb, \mathbb{P}_n(\ca)\ltimes \mathsf{H}_{n+1}(\ca)[n+2])$$ is naturaly isomorphic to the set 
$$\mathbb{R}\mbox{Der}(\cb, (M_n\circ F)^{\ast}(\mathsf{H}_{n+1}(\ca)[n+2]))$$ of derived derivations of $\cb$ with values in $(M_n\circ F)^{\ast}(\mathsf{H}_{n+1}(\ca)[n+2])$. This implies that the obstruction class $\omega_n(F)$ of $F$ corresponds to a derived derivation of $\cb$ with values in the $\cb\text{-}\cb$-bimodule $(M_n\circ F)^{\ast}(\mathsf{H}_{n+1}(\ca)[n+2])$. Moreover by the above isomorphim, the obstruction class $\omega_n(F)$ of $F$ vanishes if and only if the associated derived derivation of $\cb$ is the trivial one, see notation~\ref{hoje1}.
\end{remark}
\begin{proposition}\label{welldef}
Let $\cb$ be a cofibrant dg category. If two dg functors $F_1, F_2: \cb \longrightarrow P_n(\ca)$ become equal in the homotopy category $\mathsf{Ho}(\dgcat)(\cb, P_n(\ca))$, they give rise to isomorphic obstruction classes. In particular $\omega_n(F_1)$ vanishes if and only if $\omega_n(F_2)$ vanishes.
\end{proposition}
\begin{proof}
Notice that since every object in $\dgcat$ is fibrant (see~\ref{objfib}) and $\cb$ is cofibrant, two dg functors $F_1$ and $F_2$ become equal in $\mathsf{Ho}(\dgcat)(\cb, P_n(\ca))$ if and only if they are left homotopic. We can then construct the following diagram
$$
\xymatrix{
*+<1pc>{\cb} \ar@{>->}[d]^{\sim}_{i_0} \ar[dr]^{F_1} & & \\
I(\cb) \ar[r]^-H & P_n(\ca) \ar[r] & \mathbb{P}_n(\ca)\\
*+<1pc>{\cb} \ar@{>->}[u]^{i_1}_{\sim} \ar[ru]_{F_2} & & \,,
}
$$
where $I(\cb)$ is a cylinder object for $\cb$ and $i_0$ and $i_1$ are quasi-equivalences. Observe that the previous diagram gives rise to a zig-zag of weak equivalences in $\dgcat\downarrow \mathbb{P}_n(\ca)$ between $\omega_n(F_1)$ and $\omega_n(F_2)$, which implies that the obstruction classes are isomorphic. In particular $\omega_n(F_1)$ vanishes if and only if so does $\omega_n(F_2)$.
\end{proof}
Let us return to our `rigidification' problem: let $\ca$ be a positively graded dg category and $F:\cb \longrightarrow \mathsf{H}_0(\ca)$ a dg functor, with $\cb$ a cofibrant dg category. Consider the diagram
$$
\xymatrix{
\vdots \ar[d] &  & \vdots \ar@{->>}[d] & & \\
P_2(\ca) \ar[d] \ar@{->>}[rr]^{\sim}_{M_2} & & \mathbb{P}_2(\ca) \ar@{->>}[d] && \\
P_1(\ca) \ar[d] \ar@{->>}[rr]^{\sim}_{M_1} & & \mathbb{P}_1(\ca) \ar@{->>}[d] & &\\
P_0(\ca) \ar@{->>}[rr]^-{\sim}_-{M_0} & & \mathbb{P}_0(\ca)=\mathsf{H}_0(\ca)  && \cb \ar[ll]^-{F_0} \ar@{-->}[llu]_{F_1} \ar@{-->}[lluu]_{F_2}\,,\\
}
$$
where the left (resp. right) column is the Big (resp. small) Postnikov model for $\ca$ and the morphism between the two is the one of remark~\ref{big->small}.

Our strategy will be to try to lift $F_0: \cb \longrightarrow \mathsf{H}_0(\ca)$ to dg functors $F_n:\cb \longrightarrow \mathbb{P}_n(\ca)$ for $n=1,2,\ldots$ in sucession. If we are able to find all these lifts, there will be no difficulty in constructing the desired lift
$$F= \underset{n}{\mbox{lim}}\, F_n: \cb \longrightarrow \ca \simeq \underset{n}{\mbox{lim}}\, \mathbb{P}_n(\ca)\,.$$

For the inductive step, we have a commutative (solid) diagram as follows ($n\geq 0$)
$$
\xymatrix{
P_{n+1}(\ca) \ar[d] \ar@{->>}[rr]^{\sim}_{M_{n+1}} & & \mathbb{P}_{n+1}(\ca) \ar@{->>}[d] & \\
P_n(\ca) \ar@{->>}[rr]_{M_n}^{\sim} & & \mathbb{P}_n(\ca) & \\
& & & \cb \ar[ul]^-{F_n} \ar@{-->}@/^-1pc/[uul]_-{F_{n+1}} \ar@{-->}@/^1pc/[lllu]^-{\widetilde{F_n}}\,.
}
$$
Since $\cb$ is cofibrant and $M_n$ is a trivial fibration, there exits a lift $\widetilde{F_n}$ of $F_n$ such that $M_n \circ \widetilde{F_n}=F_n$. Moreover since $M_n$ is a quasi-equivalence, any two such lifts become equal in $\mathsf{Ho}(\dgcat)(\cb, P_n(\ca))$ and so by proposition~\ref{welldef} they give rise to isomorphic obstruction classes. In what follows, we denote by $\omega_n(F_n)$ the obstruction class of $\widetilde{F_n}$.
\begin{proposition}\label{step}
A lift $F_{n+1}$ of $F_n$, making the diagram
$$ 
\xymatrix{
\mathbb{P}_{n+1}(\ca) \ar@{->>}[d] & \\
\mathbb{P}_n(\ca) & \cb \ar[l]^-{F_n} \ar[ul]_-{F_{n+1}}
}
$$
commute, exists if and only if the obstruction class $\omega_n(F_n)$ vanishes (see~\ref{obstruction}).
\end{proposition}
\begin{proof}
Let us suppose first that $\omega_n(F_n)$ vanishes. Recall from theorem~\ref{mainthm}, that we have an homotopy fiber sequence $$ P_{n+1}(\ca) \longrightarrow P_n(\ca) \stackrel{\gamma_n}{\longrightarrow} \mathbb{P}_n(\ca)\ltimes \mathsf{H}_{n+1}(\ca)[n+2]$$ in $\mathsf{Ho}(\dgcat \downarrow \mathbb{P}_n(\ca))$. By hypothesis, the obstruction class $\omega_n(F_n)$ vanishes and so the choice of a homotopy in $\dgcat \downarrow \mathbb{P}_n(\ca)$ between $\gamma_n \circ \widetilde{F_n}$ and
$$\cb \longrightarrow \mathbb{P}_n(\ca) \longrightarrow \mathbb{P}_n(\ca)\ltimes \mathsf{H}_{n+1}(\ca)[n+2]\,\,\,\,\,\,\,\,\,(\mbox{see}~\ref{natpoint1})$$
induces a morphism in $\mathsf{Ho}(\dgcat \downarrow \mathbb{P}_n(\ca))(\cb, P_{n+1}(\ca))$. Since $\cb$ is cofibrant (and $P_{n+1}(\ca)$ is fibrant) in $\dgcat \downarrow \mathbb{P}_n(\ca)$ (see \ref{rks}), we can represent this morphism by a dg functor $\psi: \cb \longrightarrow P_{n+1}(\ca)$. Moreover, by lemma~\ref{rks1}, any two such representatives become equal in $\mathsf{Ho}(\dgcat)(\cb, P_{n+1}(\ca))$. This implies that $F_n$ and the composition
$$\xymatrix{ \cb \ar[rr]^-{M_{n+1}\circ \psi} & & \mathbb{P}_{n+1}(\ca)\ar@{->>}[r] & \mathbb{P}_n(\ca)}$$
becomes equal in $\mathsf{Ho}(\dgcat)(\cb, \mathbb{P}_n(\ca))$. Finally, by lemma~\ref{strictification}, we conclude that there exists a desired lift $F_{n+1}$ as in the proposition.

Let us now prove the converse. Suppose we have a lift $F_{n+1}$ of $F_n$ as in the proposition. Since $\cb$ is cofibrant and $M_{n+1}$ is a trivial fibration there exists a lift $\widetilde{F_{n+1}}$ of $F_{n+1}$ such that $M_{n+1}\circ \widetilde{F_{n+1}}=F_{n+1}$. Observe that $\widetilde{F_n}$ and the composition
$$ \cb \stackrel{\widetilde{F_{n+1}}}{\longrightarrow} P_{n+1}(\ca) \longrightarrow P_n(\ca)$$
becomes equal in $\mathsf{Ho}(\dgcat)(\cb, P_n(\ca))$. This implies, by theorem~\ref{mainthm} and proposition~\ref{welldef}, that the obstruction class $\omega_n(F_n)$ vanishes.
\end{proof}
Thus if it happens that at each stage of the inductive process of constructing lifts $F_n: \cb \longrightarrow \mathbb{P}_n(\ca)$, the obstruction class $\omega_n(F_n)$ vanishes, then the `rigidification' problem has a solution.
\begin{theorem}\label{final}
Let $\ca$ be a positively graded dg category and $F_0: \cb \longrightarrow \mathsf{H}_0(\ca)$ a dg functor, with $\cb$ a cofibrant dg category. If the family $\{ \omega_n(\widetilde{F_n}) \}_{n \geq 0}$ of obstruction classes vanishes, then there exists a lift $F: \cb \longrightarrow \ca$ of $F_0$, making the diagram
$$
\xymatrix{
& \ca \ar[d]^{\tau_{\leq 0}} \\
\cb \ar[r]_-{F_0} \ar[ur]^F & \mathsf{H}_0(\ca)
}
$$
commute.
\end{theorem}

\end{document}